\documentclass[12pt]{article}

\usepackage{amsmath}
\usepackage{float}
\usepackage{amsfonts}
\usepackage{amstext}
\usepackage{graphics}
\usepackage{epsf}
\usepackage[dvips]{epsfig}
\usepackage{fullpage}

\newtheorem{theorem}{Theorem}
\newtheorem{lemma}[theorem]{Lemma}
\newtheorem{corollary}[theorem]{Corollary}

\newfont{\cyr}{wncyr10 at 12.0pt}

\def \prend{\vrule depth-1pt height7pt width6pt}

\def \endpf{{\ \ \prend \medbreak}}

\title{Simultaneous Avoidance of Large Squares and Fractional Powers in Infinite Binary Words}

\author{Jeffrey Shallit\\
School of Computer Science \\
University of Waterloo \\
Waterloo, ON, N2L 3G1 \\
CANADA \\
{\tt shallit@graceland.uwaterloo.ca}
}

\begin{document}
\date{\today}
\maketitle

\begin{abstract}
In 1976, Dekking showed that there exists an infinite binary word
that contains neither squares $yy$ with $|y| \geq 4$ nor cubes $xxx$.
We show that `cube' can be replaced by any fractional power $> 5/2$.  We also
consider the analogous problem where `$4$' is replaced by any integer.
This results in an interesting and subtle hierarchy.
\end{abstract}

\section{Introduction}

       A {\it square} is a nonempty word of the form $yy$, as in the
English word {\tt murmur}.  It is easy to see that every word of 
length $\geq 4$ constructed from the symbols $0$ and $1$ contains  
a square, so it is impossible to avoid squares in infinite binary words.
However, in 1974, Entringer, Jackson, and Schatz 
\cite{Entringer&Jackson&Schatz:1974} proved the surprising fact that
there exists an infinite binary word containing no squares $yy$ with
$|y| \geq 3$.  Further, the bound $3$ is best possible.

      A {\it cube} is a nonempty word of the
form $xxx$, as in the English sort-of-word {\tt shshsh}.
An {\it overlap} is a word of the form $axaxa$, where $a$ is a single letter
and $x$ is a (possibly empty) word, as in the French word {\tt entente}.
Dekking \cite{Dekking:1976} showed that there exists an infinite
binary word that contains neither squares
$yy$ with $|y| \geq 4$ nor cubes $xxx$.  Furthermore, the bound $4$ is best possible.
He also proved that every overlap-free word contains arbitrarily large
squares.

      These two results suggest the following natural question:  for each
length $l \geq 1$, determine the fractional exponent $p$ (if it exists) such that
\begin{itemize}
\item[(a)] there is no infinite binary word simultaneously avoiding
squares $yy$ with $|y| \geq l$ and fractional powers $x^e$ with $e \geq p$;

\item[(b)] there is an infinite binary word simultaneously avoiding
squares $yy$ with $|y| \geq l$ and fractional powers $x^e$ with $e > p$?
\end{itemize}

      Here we say a word $w$ is an $e$'th power ($e$ rational)
if there exist words $y, y'\in \Sigma^*$
such that $w = y^n y'$, and $y'$ is a prefix of $y$ with
$n + |y'|/|y| = e$.  For example, the English word
{\tt abracadabra} is an ${{11} \over 7}$-power.  We say a word {\it avoids $p$ powers}
if it contains no subword of the form $y^e$ with $e \geq p$.  We say a word
{\it avoids $p^+$ powers} if it contains no subword of the form
$y^e$ with $e > p$.

      In this paper we completely resolve this question.  It turns out there
is a rather subtle hierarchy depending on $l$.  The results are summarized
in Table 1.

\begin{figure}[H]
\begin{center}
\begin{tabular}{|c|c|c|}
\hline
minimum length $l$& avoidable & unavoidable \\
of square avoided & power & power \\
\hline
$2$ & none & all \\
\hline
$3$ & $3^+$ & $3$ \\
\hline
$4,5,6$ & $(5/2)^+$ & $5/2$ \\
\hline
$ \geq 7$ & $(7/3)^+$ & $7/3$ \\
\hline
\end{tabular}
\end{center}
\caption{Summary of Results}
\end{figure}

      More precisely, we have 

\begin{theorem}

\item[(a)]  There are no infinite binary words that avoid
all squares $yy$ with $|y| \geq 2$.  

\item[(b)]  There are no infinite binary words that simultaneously
avoid all squares $yy$ with $|y| \geq 3$ and cubes $xxx$.

\item[(c)]  There is an infinite binary word that simultaneously avoids
all squares $yy$ with $|y| \geq 3$ and $3^+$ powers.

\item[(d)]  There is an infinite binary word that simultaneously
avoids all squares $yy$ with $|y| \geq 4$ and ${5 \over 2}^+$ powers.

\item[(e)]  There are no infinite binary words that simultaneously
avoid all squares $yy$ with $|y| \geq 6$ and ${5 \over 2}$ powers.

\item[(f)]  There is an infinite binary word that simultaneously
avoids all squares $yy$ with $|y| \geq 7$ and ${7 \over 3}^+$ powers.

\item[(g)]  For all $t \geq 1$, there are no infinite binary words
that simultaneously avoid all squares $yy$ with $|y| \geq t$ 
and ${7 \over 3}$ powers.  

\end{theorem}

     The result (a) is originally due to Entringer, Jackson, and Schatz
\cite{Entringer&Jackson&Schatz:1974}.  The result (b) is due to Dekking
\cite{Dekking:1976}.   The result (g) appears in a recent paper
of the author and J. Karhum\"aki \cite{Karhumaki&Shallit:2003}.
We mention them for completeness.  The remaining
results are new.

\section{Proofs of the negative results}

     We say a word avoids $(l,p)$ if it simultaneously
avoids squares $yy$ with $|y| \geq l$ and $p$ powers.

     The negative results (a), (b), and (e) can be proved
purely mechanically.  The idea is as follows.  Given $l$ and $p$, we can
create a tree $T = T(l,p)$ of all binary words avoiding $(l,p)$
as follows:  the root of $T$ is labeled $\epsilon$. If
a node is labeled $w$ and avoids $(l,p)$, then it is an internal node with two children,
where the left child is labeled $w0$ and the right child is labeled $w1$.
If it does not avoid $(l,p)$, then it is an external node (or ``leaf'').

     It is now easy to see that no infinite word avoiding $(l,p)$ exists if and only
if $T(l,p)$ is finite.  In this case, a breadth-first search will suffice to
resolve the question.  Furthermore, certain parameters of $T(l,p)$ correspond
to information about the finite words avoiding $(l,p)$:

\begin{itemize} 

\item the number of leaves $n$ is one more than the number of internal nodes,
and so $n-1$ represents the total number of finite words avoiding $(l,p)$;

\item if the height of the tree (i.e., the length of the longest path from the root
to a leaf) is $h$, then $h$ is the smallest integer such that there are no
words of length $\geq h$ avoiding $(l,p)$;

\item the internal nodes at depth $h-1$ gives the all words of maximal length
avoiding $(l,p)$.

\end{itemize}

     The following table lists $(l, p, n, h, t, S)$, where
\begin{itemize}
\item $l = |y|$, where one is trying avoiding $yy$;

\item $p$, the fractional exponent one is trying to avoid;

\item $n$, the number of leaves of $T(l,p)$;

\item $h$, the height of the tree $T(l,p)$.  

\item $t$, the number of internal nodes at depth $h-1$ in the tree.  

\item $S$, the set of labels of the internal nodes at depth $h-1$ that start
with $0$.  (The other words can be obtained simply by interchanging $0$ and $1$.)
\end{itemize}
For completeness, we give the results for the optimal exponents for
$2 \leq l \leq 7$.  As mentioned above, the case $l = 2$ is due to Entringer,
Jackson, and Schatz \cite{Entringer&Jackson&Schatz:1974} and the case
$l = 3$ is due to Dekking \cite{Dekking:1976}.

\begin{figure}[H]
\begin{center}
\begin{tabular}{|r|r|r|r|r|l|}
\hline
$l$ & $p$ & $n$ & $h$ & $t$ & $S$ \\
\hline
2 & $\infty$ & 478 & 19 & 2 & $\lbrace 010011000111001101 \rbrace$ \\
\hline
3 & 3 & 578 & 30 & 2 & $\lbrace 00110010100110101100101001100 \rbrace$ \\
\hline
4 & $5/2$ & 6860 & 84 & 4 & $\scriptscriptstyle\lbrace${\tiny 00101101001011001001101100101101001101100100110100101100100110110010110100110110011}, \\
& & & & & 
{\tiny 00110010011010010110010011011001011010011011001001101001011001001101100101101001011}$\scriptscriptstyle\rbrace$ \\
\hline
5 & $5/2$ & 15940 & 93 & 2 & $\scriptscriptstyle\lbrace${\tiny 00100101100110100101100100110110010110100110110010011010010110010011011001011010011001011011}$\scriptscriptstyle\rbrace $ \\
\hline
6 & $5/2$ & 15940 & 93 & 2 & $\scriptscriptstyle\lbrace${\tiny 00100101100110100101100100110110010110100110110010011010010110010011011001011010011001011011}$\scriptscriptstyle\rbrace$ \\
\hline
7 & $7/3$ & 3548 & 43 & 2 & $\lbrace 001011001011010011001011001101001011001011 \rbrace$ \\
\hline
\end{tabular}
\caption{Proofs of the negative results (a), (b), (e)}
\end{center}
\end{figure}

\section{Proof of (c)}
\label{proofc}

     In this section we prove that there is an infinite binary word
that simultaneously avoids $yy$ with $|y| \geq 3$ and $3^+$ powers.

     We introduce the following notation for alphabets:  $\Sigma_k := \lbrace 0, 1, \ldots, k-1 \rbrace$.

     Let the morphism $f: \Sigma_3^* \rightarrow \Sigma_2^*$ be defined as follows.
\begin{eqnarray*}
0 &\rightarrow& 0010111010 \\
1 &\rightarrow& 0010101110 \\
2 &\rightarrow& 0011101010 
\end{eqnarray*}

     We will prove

\begin{theorem}
     If $w$ is any squarefree word over $\Sigma_3$, then $f(w)$
avoids $yy$ with $|y| \geq 3$ and $3^+$ powers.
\label{square3}
\end{theorem}

\begin{proof}
We argue by contradiction.  Let $w = a_1 a_2 \cdots a_n$ be a
squarefree string such that $f(w)$ contains a square $yy$
with $|y| \geq 3$,
i.e., $f(w) = xyyz$ for some $x, z \in \Sigma_2^*$,
$y \in \lbrace \Sigma_2^{\geq 3}$.
Without loss of generality, assume that $w$ is a shortest such
string, so that $0 \leq |x|, |z| < 20$.  

Case 1:  $|y| \leq 20$.    In this case we can take $|w| \leq 5$.
To verify that $f(w)$ has no squares $yy$ with $|y| \geq 3$,
it therefore suffices to check each of the 30 possible words $w \in
\Sigma_2^5$.

Case 2: $|y| > 20$.  First, we establish the following result.

\begin{lemma}
\begin{itemize}
\item[(a)] (inclusion property)
Suppose $f(ab) = t f(c) u$ for some
letters $a, b, c \in \Sigma_2$
and strings $t, u \in \Sigma_2^*$.
Then this inclusion is trivial (that is,
$t = \epsilon$ or $u = \epsilon$).

\item[(b)] (interchange property)
Suppose there exist letters $a, b, c$ and
strings $s, t, u, v$ such that $f(a) = st$, $f(b) = uv$,
and $f(c) = sv$.  Then either $a = c$ or $b = c$.
\end{itemize}
\label{ming}
\end{lemma}

\begin{proof}
\begin{itemize}
\item[(a)]
A short computation verifies there are
no $a, b, c$ for which the equality $f(ab) = t f(c) u$ 
holds nontrivially.

\item[(b)]  This can also be verified with a short computation.
If $|s| \geq 6$, then no two distinct
letters share a prefix of length $6$.
If $|s|\leq 5$, then $|t| \geq 5$, and no two distinct letters 
share a suffix of length $5$.
\end{itemize}
\endpf
\end{proof}

     Once Lemma~\ref{ming} is established, the rest of the argument
is fairly standard.  It can be found, for example, in
\cite{Karhumaki&Shallit:2003}, but for completeness we repeat it here.

     For $i = 1, 2, \ldots, n$ define $A_i = f(a_i)$.  
Then if $f(w) = xyyz$, we can write
$$f(w) = A_1 A_2 \cdots A_n = A'_1 A''_1 A_2 \cdots A_{j-1} 
A'_j A''_j A_{j+1} \cdots A_{n-1} A'_n A''_n$$ 
where
\begin{eqnarray*}
A_1 &=& A'_1 A''_1 \\
A_j &=& A'_j A''_j \\
A_n &=& A'_n A''_n \\
x &=& A'_1 \\
y &=& A''_1 A_2 \cdots A_{j-1} A'_j = A''_j A_{j+1} \cdots A_{n-1} A'_n \\
z &=& A''_n, \\
\end{eqnarray*}
where $|A''_1|, |A''_j| > 0$.
See Figure~\ref{fig1}.

\begin{figure}[H]
\begin{center}
\input cube1.pstex_t
\end{center}
\caption{The string $xyyz$ within $f(w)$ \protect\label{fig1}}
\end{figure}

     If $|A''_1| > |A''_j|$, then $A_{j+1} = f(a_{j+1})$ is a subword
of $A''_1 A_2$, hence a subword of $A_1 A_2 = f(a_1 a_2)$.  Thus
we can write $A_{j+2} = A'_{j+2} A''_{j+2}$ with
$$ A''_1 A_2 = A''_j A_{j+1} A'_{j+2}.$$ 
See Figure~\ref{fig2}.

\begin{figure}[H]
\begin{center}
\input cube2.pstex_t
\end{center}
\caption{The case $|A''_1| > |A''_j|$ \protect\label{fig2}}
\end{figure}

But then,
by Lemma~\ref{ming} (a), either $|A''_j| = 0$,
or $|A''_1| = |A''_j|$, or $A'_{j+2}$ is a not
a prefix of any $f(d)$.  All three conclusions are impossible.

     If $|A''_1| < |A''_j|$, then $A_2 = f(a_2)$ is a subword of
$A''_j A_{j+1}$, hence a subword of $A_j A_{j+1} = f(a_j a_{j+1})$.
Thus we can write $A_3 = A'_3 A''_3$ with
$$ A''_1 A_2 A'_3 = A''_j A_{j+1} .$$  
See Figure~\ref{fig3}.

\begin{figure}[H]
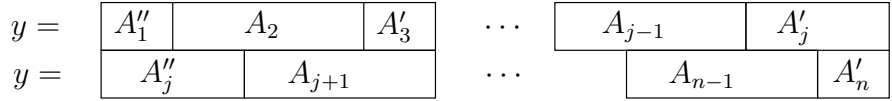

\begin{center}
\input cube3.pstex_t
\end{center}
\caption{The case $|A''_1| < |A''_j|$ \protect\label{fig3}}
\end{figure}

By Lemma~\ref{ming} (a), either $|A''_1| = 0$ or $|A''_1| = |A''_j|$
or $A'_3$ is not a prefix of any $f(d)$.  Again, all three conclusions
are impossible.

     Therefore $|A''_1| = |A''_j|$.  
Hence $A''_1 = A''_j$, $A_2 = A_{j+1}$, $\ldots$, $A_{j-1} = A_{n-1}$,
and $A'_j = A'_n$.  Since $h$ is injective, we have
$a_2 = a_{j+1}, \ldots, a_{j-1} = a_{n-1}$.
It also follows that $|y|$ is divisible by $10$ and
$A_j = A'_j A''_j = A'_n A''_1$.   But by Lemma~\ref{ming} (b), either
(1) $a_j = a_n$ or (2) $a_j = a_1$.  In the first case,
$a_2 \cdots a_{j-1} a_j = a_{j+1} \cdots a_{n-1} a_n$, so
$w$ contains the square $(a_2 \cdots a_{j-1} a_j)^2$, a contradiction.  In the
second case, $a_1 \cdots a_{j-1} = a_j a_{j+1} \cdots a_{n-1}$, so
$w$ contains the square $(a_1 \cdots a_{j-1})^2$, a contradiction.

     It now follows that if $w$ is squarefree then $f(w)$ avoids squares
$yy$ with $|y| \geq 3$.   

     It remains to see that $f(w)$ avoids $3^+$ powers.   If $f(w)$ contained
$x^e$ for some fractional exponent $e > 3$, then it would contain $x^2$, so
from above we have $|x| \leq 2$.  Thus it suffices to show that $f(w)$ avoids
the words $0000, 1111, 0101010, 1010101$.  This can be done by a short
computation.
\endpf
\end{proof}

\begin{corollary}
    There is an infinite binary word avoiding squares $yy$ with $|y| \geq 3$
and $3^+$ powers.
\end{corollary}

\begin{proof}
      As is very well-known, there are infinite squarefree words over $\Sigma_3$
\cite{Thue:1906,Berstel:1995}.
Take any such word $\bf w$ (for example, the fixed point of the morphism
$2 \rightarrow 210$, $1 \rightarrow 20$, $0 \rightarrow 1$), and apply the
map $f$.  The resulting word $f({\bf w})$ avoids $(3, 3^+)$.
\endpf
\end{proof}

     It may be of some interest to explain how the morphism $f$ was discovered.
We iteratively generated all words of length
$1, 2, 3, \ldots$ (up to some bound) that avoid $(3, 3^+)$.    We then
guessed such words were the image of a $k$-uniform morphism applied to a
squarefree word over $\Sigma_3$.  For values of $k = 2,3, \ldots$, we broke up
each word into contiguous blocks of size $k$, and discarded any word for which
there were more than $3$ blocks.  For certain values of $k$, this procedure
eventually resulted in $0$ words fitting the criteria.  At this point we knew
a $k$-uniform morphism cannot work, so we increased $k$ and started over.  Eventually
a $k$ was found for which the number of such words appeared to increase without
bound.  We then examined the possible sets of $3$ $k$-blocks to see if any
satisfied the requirements of Lemma~\ref{ming}.  This gave our candidate morphism
$f$.

\begin{theorem}
     Let $A_n$ denote the number of binary words of length $n$ avoiding $yy$ with $|y| \geq 3$
and $3^+$ powers.  Then $A_n = \Omega(1.01^n)$ and
$A_n = O(1.49^n)$.
\label{avoid3}
\end{theorem}

\begin{proof}
Grimm \cite{Grimm:2001}
has shown there are $\Omega(\lambda^n)$ squarefree words over $\Sigma_3$, where
$\lambda = 1.109999$.  Since the map $f$ is $10$-uniform,
it follows that $A_n = \Omega(\lambda^{n/10}) = \Omega(1.01^n)$.

For the upper bound, we reason as follows.  The set of binary words of length $n$ 
avoiding $yy$ with $|y| \geq 3$ and $3^+$ powers is a subset of the set of binary
words avoiding $0000$ and $1111$.  The number $A'_n$ of words avoiding
$0000$ and $1111$ satisfies the linear recurrence $A'_n = A'_{n-1} + A'_{n-2} + A'_{n-3}$
for $n \geq 4$.  It follows that $A'_n = O(\alpha^n)$, where $\alpha$ is the largest
zero of $x^3 - x^2 -x -1$, the characteristic polynomial of the recurrence.
Here $\alpha < 1.84$, so $A_n = O(1.84^n)$.

This reasoning can be extended using a symbolic algebra package such as Maple.
Noonan and Zeilberger \cite{Noonan&Zeilberger:1999}
have written a Maple package {\tt DAVID\_IAN} that allows one to specify a list $L$
of forbidden words, and computes the generating function enumerating words avoiding
members of $L$.  We used this package for a list $L$ of $62$ words of length $\leq 12$:
$$ 0000, 1111, \ldots,  111010111010 $$
obtaining a characteristic polynomial of degree $67$ with dominant zero $\doteq 1.4895$.
\endpf
\end{proof}

\section{Proof of (e)}
\label{proofe}

     In this section we prove that there is an infinite binary word
that simultaneously avoids $yy$ with $|y| \geq 4$ and ${5 \over 2}^+$ powers.

     Let $g_1: \Sigma_8^* \rightarrow \Sigma_2^*$ be defined as follows.
\begin{eqnarray*}
0 &\rightarrow& 0011010010110 \\
1 &\rightarrow& 0011010110010 \\
2 &\rightarrow& 0011011001011 \\
3 &\rightarrow& 0100110110010 \\
4 &\rightarrow& 0110100101100 \\
5 &\rightarrow& 1001101011001 \\
6 &\rightarrow& 1001101100101 \\
7 &\rightarrow& 1010011011001
\end{eqnarray*}

     Let $g_2: \Sigma_4^* \rightarrow \Sigma_8^*$ be defined as follows.
\begin{eqnarray*}
0 &\rightarrow& 03523503523453461467 \\
1 &\rightarrow& 03523503523453467167 \\
2 &\rightarrow& 16703523503523461467 \\
3 &\rightarrow& 03523503523461467167 
\end{eqnarray*}

     Let $g_3: \Sigma_3^* \rightarrow \Sigma_4^*$ be defined as follows.
\begin{eqnarray*}
0 &\rightarrow& 010203 \\
1 &\rightarrow& 010313 \\
2 &\rightarrow& 021013
\end{eqnarray*}

     Finally, define $g: \Sigma_3^* \rightarrow \Sigma_2^*$ by
$g = g_1 \circ g_2 \circ g_3$.  Note that $g$ is $1560$-uniform.

     We will prove

\begin{theorem}
     If $w$ is any squarefree word over $\Sigma_3$, then $g(w)$
avoids $yy$ with $|y| \geq 4$ and ${5 \over 2}^+$ powers.
\end{theorem}

\begin{proof}
     The proof is very similar to the proof of Theorem~\ref{square3}, and
we indicate only what must be changed.

      First, it can be checked that
Lemma~\ref{ming} also holds for the morphism $g$.

      As before, we break the proof up into two parts:  the case where
$g(w) = xyyz$ for some $y$ with $4 \leq |y| \leq 2 \cdot 1560$, and the
case where $g(w) = xyyz$ for some $y$ with $|y| \geq 2 \cdot 1560$.
The former can be checked by examining the image of the 30 squarefree words
in $\Sigma_3^5$ under $g$.  The latter is handled as we did in the proof
of Theorem~\ref{square3}.  We checked these conditions with programs written in
Pascal; these are available from the author on request.
\endpf
\end{proof}

\begin{corollary}
    There is an infinite binary word avoiding squares $yy$ with $|y| \geq 4$
and ${5 \over 2}^+$ powers.
\end{corollary}

     It may be of some interest to explain how the morphisms $g_1$, $g_2$, $g_3$,
were discovered.

     We used a procedure analogous to that described above in
Section~\ref{proofc}.   However, since it was not feasible to generate
all words avoiding $(4, {5\over2}^+)$ and having at most $3$ contiguous
blocks of length $1560$, we increased the alphabet size and and tried
various $k$-blocks until we found a combination of alphabet size and block
size for which the number of words appeared to increase without bound.
We then obtained a number of possible candidates for blocks.

     Next, we determined the necessary avoidance properties of the blocks
given by images of letters under $g_1$.  For example, $g_1(0)$ cannot be
followed by $g_1(1)$, because this results in the subword $000$, which
is a 3rd power (and $3 > 2.5$).  The blocks that must be avoided include
all words with squares, and
$$
01,02,04,05,06,07, 10,12,13,17,
20,21,24,25,26,27, 30,31,32,36,37,
40,41,42,43,47, $$
$$ 51,54,56,57,
60,62,63,64,65, 72,73,74,75,76,
034,145,153,161,353,450,452,535,615,616,
$$
$$
714,715,
2346703,5234670,5234671,
53467035,
6703523461,
2346146703503, 5234614670350
$$
This list was computed purely mechanically, and
it is certainly possible that this list is not exhaustive.

     We now iterated our guessing procedure, looking for a candidate
uniform morphism that creates squarefree words avoiding the patterns in the list
above.    This resulted in the $20$-uniform morphism $g_2$.

      We then computed the blocks that must be avoided for $g_2$.  This was
done purely mechanically.  Our procedure suggested that arbitrarily large
blocks must be avoided, but luckily they (apparently) had a simple finite
description:  namely, we must avoid $12$, $23$, $32$, and blocks of the
form $2x0x1$ and $3x1x0$ for all nonempty words $x$, in addition to words
with squares.  

     We then iterated our guessing procedure one more time, looking for a candidate
uniform morphism that avoids {\it these} patterns.  This gave us the morphism
$g_3$.  

     Of course, once the morphism $g = g_1 \circ g_2 \circ g_3$ is discovered,
we need not rely on the list of avoidable blocks; we can take the morphism
as given and simply verify the properties of inclusion and interchange
as in Lemma~\ref{ming}.

\begin{theorem}
     Let $B_n$ denote the number of binary words of length $n$ avoiding $yy$ with $|y| \geq 4$
and ${5 \over 2}^+$ powers.  Then $B_n = \Omega(1.000066^n)$ and
$B_n = O(1.122^n)$.
\label{avoid4}
\end{theorem}

\begin{proof}
      The proof is analogous to that of Theorem~\ref{avoid3}.  We use the fact
that $g$ is $1560$-uniform, which, combined with the result of 
Grimm \cite{Grimm:2001},
gives the bound $1.109999^{1/1560} \doteq 1.000066899$.

      For the upper bound, we again use the Noonan-Zeilberger Maple package.
We used the $54$ patterns corresponding to words of length $\leq 20$.
This gave us a polynomial of degree $27$ with dominant zero
$\doteq 1.12123967$.
\endpf
\end{proof}

\section{Proof of (f)}     
\label{proofg}

     In this section we prove that there is an infinite binary word
that simultaneously avoids $yy$ with $|y| \geq 7$ and ${7 \over 3}^+$ powers.

     Let $h_1: \Sigma_5^* \rightarrow \Sigma_2^*$ be defined as follows.
\begin{eqnarray*}
0 &\rightarrow& 00110100101100 \\
1 &\rightarrow& 00110100110010 \\
2 &\rightarrow& 01001100101100 \\
3 &\rightarrow& 10011011001011 \\
4 &\rightarrow& 11010011011001
\end{eqnarray*}

     Let $h_2: \Sigma_3^* \rightarrow \Sigma_5^*$ be defined as follows.
\begin{eqnarray*}
0 &\rightarrow& 032303241403240314  \\
1 &\rightarrow& 032314041403240314 \\
2 &\rightarrow& 032414032303240314 
\end{eqnarray*}

     Finally, define $h: \Sigma_3^* \rightarrow \Sigma_2^*$ by
$h = h_1 \circ h_2$.  Note that $h$ is $252$-uniform.

     We will prove

\begin{theorem}
     If $w$ is any squarefree word over $\Sigma_3$, then $h(w)$
avoids $yy$ with $|y| \geq 7$ and ${7 \over 3}^+$ powers.
\end{theorem}

\begin{proof}
     Again, the proof is quite similar to that of Theorem~\ref{square3}.
We leave it to the reader to verify that the inclusion and interchange
properties hold for $h$, and that the image of all the squarefree words
of length $\leq 5$ are free of squares $yy$ with $|y| < 7$ and
${7 \over 3}^+$ powers.
\endpf
\end{proof}

\begin{corollary}
    There is an infinite binary word avoiding squares $yy$ with $|y| \geq 7$
and ${7 \over 3}^+$ powers.
\end{corollary}

     The morphisms $h_1, h_2$ were discovered using the heuristic procedure
mentioned in Section~\ref{proofc}.  The avoiding blocks for $h_1$ were
heuristically discovered to include
$$ 01, 02, 10, 12, 13, 20, 21, 34, 42, 43, 304, 23031, 24041, 231403141,  232403241$$
as well as blocks containing any squares.
Then $h_2$ was constructed to avoid these blocks.

\begin{theorem}
     Let $C_n$ denote the number of binary words of length $n$ avoiding $yy$ with $|y| \geq 7$
and ${7 \over 3}^+$ powers.  Then $C_n = \Omega(1.0004^n)$ and
$C_n = O(1.162^n)$.
\end{theorem}

\begin{proof}
     The proof is very similar to that of Theorems~\ref{avoid3} and \ref{avoid4}.

     For the lower bound, note that $h$ is $252$-uniform.  This, combined with
the bound of Grimm \cite{Grimm:2001}, gives a lower bound of $\Omega(\lambda^n)$ for
all $\lambda < 1.109999^{1/252} \doteq 1.0004142$.

     For the upper bound, we again used the Noonan-Zeilberger Maple package.
We avoided $58$ words of length $\leq 20$.  This resulted in a polynomial of
degree $26$, with dominant zero $\doteq 1.1615225$.
\endpf
\end{proof}

%
%
%
%
%
%
%
%
%
%
%
%
%
%
%

\section{Enumeration results}

     In this section we provide a table of the first values of the sequences
$A_n$, $B_n$, and $C_n$, defined in Sections~\ref{proofc}, \ref{proofe}, and
\ref{proofg}, for $1 \leq n \leq 25$.

\begin{figure}[H]
\begin{center}
\begin{tabular}{|r|r|r|r|}
\hline
$n$ & $A_n$ & $B_n$ & $C_n$ \\
\hline
0 & 1 & 1 & 1 \\
1 & 2 & 2 & 2 \\
2 & 4 & 4 & 4 \\
3 & 8 & 6 & 6 \\
4 & 14 & 10 & 10 \\
5 & 26 & 16 & 14 \\
6 & 42 & 24 & 20 \\
7 & 68 & 36 & 30 \\
8 & 100 & 46 & 38 \\
9 & 154 & 64 & 50 \\
10 & 234 & 74 & 64 \\
11 & 356 & 88 & 86 \\
12 & 514 & 102 & 108 \\
13 & 768 & 114 & 136 \\
14 & 1108 & 124 & 164 \\
15 & 1632 & 140 & 196 \\
16 & 2348 & 160 & 226 \\
17 & 3434 & 178 & 264 \\
18 & 4972 & 198 & 322 \\
19 & 7222 & 212 & 384 \\
20 & 10356 & 230 & 436 \\
21 & 14962 & 256 & 496 \\
22 & 21630 & 294 & 578 \\
23 & 31210 & 342 & 674 \\
24 & 44846 & 366 & 754 \\
25 & 64584 & 392 & 850 \\
\hline
\end{tabular}
\end{center}
\caption{Values of $A_n$, $B_n$, $C_n$, $0 \leq n \leq 25$}
\end{figure}

\section{Acknowledgments}

     I would like to thank Jean-Paul Allouche and Matthew Nichols for having
suggested the problem, and Narad Rampersad for independently verifying the properties
of the morphism $g$.

\end{document}